\documentstyle [11pt]{article}

\newtheorem{propo}{{\bf Proposition}}[section]

\newtheorem{lemma}[propo]{{\bf Lemma}} \newtheorem{theor}[propo]{{\bf Theorem}}

\begin{document}

\centerline{\bf TWO GENERATOR SUBALGEBRAS OF LIE ALGEBRAS}

\bigskip

\centerline {Kevin Bowman} \centerline {Department of Physics,
Astronomy and Mathematics} \centerline {University of Central
Lancashire} \centerline {Preston PR1 2HE, England}
\smallskip

\centerline {David A. Towers} \centerline {Department of
Mathematics, Lancaster University} \centerline {Lancaster LA1 4YF,
England}

\centerline {and}

\centerline {Vicente R. Varea \footnote[1]{Supported by DGI Grant
BFM2000-1049-C02-01}} \centerline {Department of Mathematics,
University of Zaragoza} \centerline {Zaragoza, 50009 Spain}

\bigskip

{\sl Abstract}  In \cite{tho} Thompson showed that a finite group
$G$ is solvable if and only if every two-generated subgroup is
solvable (Corollary 2, p. 388). Recently, Grunevald et al.
\cite{gknp} have shown that the analogue holds for
finite-dimensional Lie algebras over infinite fields of
characteristic greater than $5$. It is a natural question to ask
to what extent the two-generated subalgebras determine the
structure of the algebra. It is to this question that this paper
is addressed. Here, we consider the classes of strongly-solvable
and of supersolvable Lie algebras, and the property of
triangulability.
\par


\par

AMS 2000 {\sl Mathematics subject classification} 17B05, 17B30

\section {Introduction.}
Let $\mathcal {P}$ be a certain property that a subalgebra of a
Lie algebra may possess. A task is to obtain information on the
structure of a Lie algebra $L$ all of whose two-generated {\em
proper} subalgebras possess the property $\mathcal{P}$. Given a
subalgebra $S$ of $L$, we distinguish two types of properties that
$S$ may possess: one type is that $S$ belong to a certain class
$\mathcal{C}$ of Lie algebras, in this case the class
$\mathcal{C}$ will be identified with the property; the other one
is that $S$ be immersed in $L$ in a certain way. In this paper we
will consider properties of these two types. Of the first type, we
will consider the classes of strongly-solvable and supersolvable
Lie algebras and of the other type the property that $S$ be
triangulable on $L$.

In section 2, we give some preliminary results and we collect some
known ones on important classes of Lie algebras. Most of them will
be used in the sequel

In section 3, we consider the classes of strongly-solvable and of
supersolvable Lie algebras. We prove that if $L$ is a solvable Lie
algebra all of whose two-generated proper subalgebras are strongly solvable 
(resp. supersolvable) then either $L$ is two-generated
and every proper subalgebra of $L$ is strongly solvable (resp.
supersolvable) or else $L$ is strongly solvable (resp.
supersolvable). In the strongly-solvable case, we also prove that
if $L$ is non-solvable, $F$ is infinite and $\mathrm{char}(F)>5$,
then $L$ is two-generated, every proper subalgebra of $L$ is
strongly solvable and $L/{\phi(L)}$ is simple (where ${\phi(L)}$
denotes the largest ideal of $L$ contained in every maximal
subalgebra of $L$).

In section 4, we prove that a solvable Lie algebra $L$ is
triangulable whenever every two-generated {\em proper} subalgebra
of $L$ is triangulable on $L$. Also, we prove that if every proper
subalgebra of a Lie algebra $L$ is triangulable on $L$ but $L$
itself is not, then $L$ is two-generated and $L/\phi(L)$ is
simple. Moreover, we give some information on the structure of a
simple Lie algebra $L$ all of whose proper subalgebras are
triangulable on $L$. In particular we obtain that such a Lie
algebra is two-generated.

Throughout $L$ will denote a finite-dimensional Lie algebra over a
field $F$. There will be no assumptions on $F$ other than those
specified in individual results. We shall call $L$ {\em
supersolvable} if there is a chain $0 = L_{0} \subset L_{1}
\subset \ldots \subset L_{n-1} \subset L_{n} = L$, where $L_{i}$
is an $i$-dimensional ideal of $L$. We shall call $L$ {\em
strongly solvable} if its derived subalgebra $L^2$ is nilpotent.
It is well known that $$ {\mathrm
supersolvable}\,\Longrightarrow\, {\mathrm
strongly-solvable}\,\Longrightarrow\,{\mathrm solvable} $$ For
algebraically closed fields of characteristic zero, these three
classes of Lie algebras coincide (Lie's theorem). For fields of
characteristic zero, every solvable Lie algebra is strongly solvable. 
There are well-known examples of solvable Lie algebras
over algebraically closed fields of non-zero characteristic which
are not supersolvable (see for instance \cite[p.53]{J} or
\cite{bn}). Nevertheless, if the ground field $F$ is algebraically
closed then every strongly-solvable Lie algebra is supersolvable
(see \cite[Lemma 2.4] {bn}).

A subalgebra $S$ of $L$ is said to be {\em triangulable} on $L$ if
${\rm ad}_{L} S = \{{\rm ad}_{L} x \mid x \in S\}$ is a Lie
algebra of linear transformations of $L$ which is triangulable
over the algebraic closure of $F$ (equivalently, if every element
of $S^2$ acts nilpotently on $L$, see \cite{w}). Some questions
regarding triangulabilty of linear Lie algebras have been
considered in \cite{gu}.

The symbol $\dot{+}$ will denote a vector space direct sum. We say
that $L$ is {\em almost abelian} if $L = L^{2}\dot{+} Fx$ with
${\rm ad}\,x$ acting as the identity map on the abelian ideal
$L^{2}$; $L$ is {\em quasi-abelian} if it is abelian or almost
abelian. Quasi-abelian Lie algebras are precisely those in which
every subspace is a subalgebra.

\section{Preliminary Results}

The {\em Frattini subalgebra} of $L$, denoted by $F(L)$, is the
intersection of all maximal subalgebras of $L$. It is known that
$F(L)$ need not be an ideal of $L$, even for algebraically closed
fields (see \cite{v3}); the {\em Frattini ideal} of $L$, denoted
by $\phi (L)$, is the largest ideal of $L$ contained in $F(L)$. We
say that $L$ is {\em $\phi$ - free} if $\phi (L) = 0$. Clearly,
$L/\phi(L)$ is $\phi$-free.

The following is straightforward.

\begin{lemma}\label{l:lemma2}
A Lie algebra $L$ is two-generated if and only if $L/\phi(L)$ is
two-generated.
\end{lemma}

A class $\mathcal{C}$ of Lie algebras is said to be {\em
saturated} if $L\in\mathcal{C}$ whenever
$L/\phi(L)\in\mathcal{C}$. It is well-known that the classes of
solvable, strongly-solvable, supersolvable and nilpotent Lie
algebras are saturated (see \cite{bn} and \cite{to})

For short, we will say that the property $\mathcal{P}$ satisfies
condition (*) if for every Lie algebra $L$ all of whose
two-generated {\it proper} subalgebras of $L$ possess the property
$ \mathcal{P}$, either $L$ itself possesses the property $
\mathcal{P}$ or $L$ is two-generated. Next, we collect some known
results on classes of Lie algebras which satisfy this condition.

{\bf Theorem 0} \begin{enumerate}

\item The classes of abelian, nilpotent
and quasi-abelian Lie algebras satisfy condition (*).

\item The class of simple
(including the one-dimensional) Lie algebras satisfies condition
(*). If ${\mathrm char}(F)\neq 2,3$ and every two-generated proper
subalgebra of $L$ is either simple or one-dimensional, then every
subalgebra of $L$ of dimension $>1$ is simple and $L$ is
two-generated.

\item If $F$ is infinite and ${\mathrm char}(F)>5$, then the
class of solvable Lie algebras satisfies condition (*)

\end{enumerate}

{\bf Proof.} (1): Clearly, the class of abelian Lie algebras
satisfies condition (*). Now, suppose that every two-generated
subalgebra of $L$ is nilpotent. Let $x \in L$ and choose any $y
\in L$. Then $<x, y>$ is nilpotent, so $y({\rm ad}\,x)^{n} = 0$
for some $n$. Hence $y({\rm ad}\,x)^{d} = 0$, where  $d =$ dim
$L$.  This holds for all $y \in L$, so ${\rm ad}\,x$ is nilpotent
for all $x \in L$. It follows from Engel's Theorem that $L$ is
nilpotent. To prove the last assertion in (1), suppose that every
two-generated subalgebra of $L$ is quasi-abelian. Then every
two-dimensional subspace of $L$ is a subalgebra of $L$, from which
it follows that every subspace of $L$ is a subalgebra of $L$, and
hence that $L$ is quasi-abelian.

(2): This is proved in \cite[Proposition 3.2]{v} and \cite[Theorem
4]{v2}.

(3): It is an immediate consequence of a result of Grunewald,
Kunyavskii, Nikolova and Plotkin in \cite{gknp}.

The structure of solvable minimal non-abelian Lie algebras has
been fully described by Stitzinger in \cite[Theorem 1]{stit}. If
$L$ is such a Lie algebra, then $L=A\dot{+} Fx$, where $A$ is an
abelian ideal of $L$ and ${\rm ad}\,x$ acts irreducibly on $A$.
Non-solvable minimal non-abelian Lie algebras have been studied by
Farnsteiner in \cite{F} and Gein in \cite{gmin} (see also Elduque
\cite{e}).

We finish this section by collecting some known results on the
structure of a minimal non-$\mathcal{C}$ Lie algebra for several
classes $\mathcal{C}$.

Minimal non-quasi-abelian Lie algebras over a field of
characteristic different from 2 and 3 have been studied by Gein in
\cite[Proposition 3]{gmod} and Varea in \cite[Theorem 2.2 and
Corollary 2.4] {vhirosh}. If $L$ is such a Lie algebra, then one
of the following occurs: (a) $L$ is solvable minimal non-abelian;
(b) $L\cong\rm{sl}(2)$; (c) $L$ is simple minimal non-abelian; or
(d) $L$ has a basis $a_1,a_2,x$ with product given by one of the
following rules:

\begin{enumerate}\item $[a_1,a_2]=0, [a_1,x]=a_1, [a_2,x]=\alpha a_2$, $1\neq
\alpha\in F$\item $[a_1,a_2]=0, [a_1,x]=a_1, [a_2,x]=a_1+ a_2$.

\end{enumerate}

Minimal non-solvable Lie algebras over an algebraically closed
field of prime characteristic have been studied by Varea in
\cite{vmin}

Minimal non-nilpotent Lie algebras have been studied by Stitzinger
in \cite{stit}, Gein and Kuznecov \cite{gk}, Towers \cite{tomin}
and Farnsteiner \cite{F}. In Gein \cite{Gbook}, it is proved that
if $L$ is simple and minimal non-nilpotent then the intersection
of two distinct maximal subalgebras of $L$ is zero and $L$ has no
non-zero ad-nilpotent elements. For the readers' convenience we
include here a proof of Gein's result. (Our proof is different
from the one given in \cite{Gbook}.)

\begin{theor}{\bf{(Gein \cite{Gbook})}} \label{t:G} Let $L$ be a simple Lie algebra over an arbitrary field
$F$. Assume that every proper subalgebra is nilpotent. Then the
following hold:
\begin{enumerate}
\item $M_{1} \cap M_{2} = 0$ for every pair of different maximal subalgebras
$M_{1}$ and $M_{2}$ of $L$;
\item $L$ has no non-zero ad-nilpotent elements; and
\item $L$ is two-generated.
\end{enumerate}
\end{theor}

{\em Proof.} (1): Let $M$ be a maximal subalgebra of $L$. Assume
that there exists a proper subalgebra $S$ of $L$ not contained in
$M$ such that $M \cap S \neq 0$. Choose $S$ such that $\dim M \cap
S$ is maximal. Nilpotency of $M$ implies that $N_{M}(M \cap S)
\neq M \cap S$. Nilpotency of $S$ implies that $N_{S}(M \cap S)
\neq M \cap S$. Let $T$ be the subalgebra of $L$ generated by
$N_{M}(M \cap S)$ and $N_{S}(M \cap S)$. Since $M \cap S$ is a
non-zero ideal of $T$, it follows from the simplicity of $L$ that
$T \neq L$. Moreover we have $S \cap M < T \cap M$, which
contradicts our choice of $S$.

(2): Let $0 \neq x \in L$ be ad-nilpotent. Let $H$ be a maximal
subalgebra of $L$ containing $x$. Let $L = L_0\dot{+} L_{1}$ be
the Fitting decomposition of $L$ relative to $H$. Since $H$ is a
maximal subalgebra of $L$, we have $H=L_0$. Since $x$ acts
nilpotently on $L_{1}$, there exists $0 \neq y \in L_{1}$ such
that $[x, y] = 0$. Therefore $C_{L}(x)$ is not contained in $H$.
Take a maximal subalgebra $M$ of $L$ containing $C_{L}(x)$. We
find that $H \cap M \neq 0$ and $H \neq M$, which contradicts (1).

(3): Assume $L \neq <x,y>$ for every $x$, $y \in L$. Let $M$ be a
maximal subalgebra of $L$. Take $0 \neq x \in M$ and $y\in L$,
$y\not\in M$. There exists a maximal subalgebra $S$ of $L$
containing $<x,y>$. We have $S \cap M \neq 0$, which contradicts
(1).
\medskip

{\it Conjecture:} Every simple minimal non-solvable Lie algebra is
two-generated.

From \cite{gknp} (see Theorem 0(3)) it follows that the conjecture
is true for infinite fields of characteristic greater than $5$. In
section 4, we will prove that a Lie algebra $L$ all of whose
proper subalgebras are triangulable on $L$ is two-generated.

\section{The classes of strongly-solvable and supersolvable Lie algebras}

In this section we consider the classes of strongly-solvable and
supersolvable Lie algebras. We recall that the Lie algebra with
basis $a,b,c$ and product given by $[a,b]=c$, $[a,c]=[b,c]=0$ is
called the {\em three-dimensional Heisenberg algebra}.

The structure of solvable $\phi$-free minimal
non-strongly-solvable Lie algebras was determined in \cite{bt} and
is given below.

\begin{theor}\label{t:bt} {\sl Let $L$ be a solvable $\phi$-free minimal non-strongly-solvable Lie algebra. Then $F$ has characteristic $p
> 0$ and $L = A \dot{+} B$ is a semidirect sum, where $A$ is the
unique minimal ideal of $L$, dim $A \geq 2$, $A^{2} = 0$, and
either $B = M \dot{+} Fx$, where $M$ is an abelian minimal ideal
of $B$ (type I), or $B$ is the three-dimensional Heisenberg
algebra (type II).}
\end{theor}

By using Theorem \ref{t:bt} we prove the following.

\begin{propo}\label{p:1}
Let $L$ be solvable $\phi$-free minimal non-strongly-solvable Lie
algebra. Then $L$ is two-generated.
\end{propo}

{\em Proof.} The structure of $L$ is described in Theorem
\ref{t:bt} above.

Suppose first that $L$ is of type I. Clearly, $[A,B]$ is an ideal
of $L$. So that either $[A,B]=A$ or $[A,B]=0$. In the latter case,
we have that $B$ is also an ideal of $L$, which is a
contradiction. Therefore $[A,B]=A$. On the other hand, we have
that $[x,M]$ is an ideal of $B$ contained in $M$. So, either
$[x,M]=M$ or $[x,M]=0$. In the latter case, we have that $B$ is
abelian. This yields that $L^2\leq A$ and hence $L^2$ is abelian,
a contradiction. Therefore $[x,M]=M$. Then we have that $B^2=M$
and $L^2=A+B^2=A+M$. If $A ({\rm ad}\, m)^{2} = 0$ for all $m \in
M$ it is easy to see that $({\rm ad}\, n)^{3} = 0$ for all $n \in
L^2$. But then $L^2$ is nil and hence nilpotent, by Engel's
theorem, a contradiction. Choose $m \in M$ such that $A ({\rm
ad}\, m)^{2} \neq 0$, and then $a \in A$ such that $[[a,m],m] \neq
0$. Now $B = <m, x>$. Put $D = <a+x, m>$. Then $[[a+x,m],m] =
[[a,m],m] \in D$, whence $D \cap A \neq 0$. Clearly $L = A + D$ so
$D \cap A$ is an ideal of $L$. It follows that $A \subseteq D$,
giving $D = A + <x, m> = L$ and $L$ is two-generated.

So suppose that $L$ is of type II. Then $B$ has a basis $c,s,x$
with product given by $[s,x]=c$, $[c,x]=[c,s]=0$. Put $C =
C_{A}(c)$. Then $C$ is an ideal of $L$ and so $C = A$ or $C = 0$.
The former implies that $L^{2} = A + Fc$ is abelian, a
contradiction, so $C = 0$. Let $D = <a + s, x>$. Then $[a,x] + c
\in D$. If $[a,x] = 0$, then $c \in D$. But this implies that $0
\neq [a,c] \in D \cap A$, whence $L = D$ as above. So suppose that
$[a,x] \neq 0$. This yields that $[[a,x],x] \in D$, from which
$[[a,x],x] \neq 0$ would give $L = D$ again. So we have that
$[a,x] \neq 0$ but $[[a,x],x] = 0$. Put $E = <[a,x] + s, x>$. Then
$c = [[a,x] + s,x] \in E$, whence $[[a,x],c] \in E$. But
$[[a,x],c] \neq 0$ since $[a,x] \neq 0$, giving $E \cap A \neq 0$
from which $L = E$. The proof is complete.

Now we are able to prove the following result.

\begin{theor}\label{t:theor1} Let $L$ be a Lie algebra such that every two-generated
proper subalgebra is strongly solvable.
\begin{enumerate}
\item Assume that $L$ is solvable but not strongly solvable. Then
every proper subalgebra of $L$ is strongly solvable and $L$ is
two-generated.

\item Assume that $L$ is non-solvable, $F$ is infinite and $\mathrm{char}(F)>5$.
Then $L$ is two-generated, every proper subalgebra of $L$ is
strongly solvable and $L/{\phi(L)}$ is simple.
\end{enumerate}
\end{theor}

{\em Proof.} (1): Let $S$ be a non-strongly-solvable subalgebra of
$L$ of minimal dimension. Clearly, $S$ is minimal
non-strongly-solvable. We have that $S/\phi(S)$ is $\phi$-free,
solvable and minimal non-strongly-solvable. By Proposition
\ref{p:1} it follows that $S/\phi(S)$ is two-generated, whence $S$
is also two-generated, by Lemma \ref{l:lemma2}. By our hypothesis
it follows that $S=L$. This yields that $L$ is minimal
non-strongly-solvable and two-generated. The proof of (1) is
complete.

(2): By using Theorem 0(3) we obtain that every non-solvable
subalgebra of $L$ (including $L$ itself) is two-generated. From
this and our hypothesis it follows that every proper subalgebra of
$L$ is solvable. Now let $S$ be a proper subalgebra of $L$. Assume
that $S$ is not strongly solvable. Then, by (1) it follows that
$S$ is two-generated. But then, by our hypothesis $S$ is strongly
solvable, a contradiction. Put $\bar{L}=L/\phi(L)$. We have that
$\bar{L}$ is non-solvable but every proper subalgebra of $\bar{L}$
is solvable. On the other hand, since $\bar{L}$ is $\phi$-free we
have that $\bar{L}=N+S$, $N\cap S=0$, where $N$ is the largest
nilpotent ideal of $\bar{L}$ and $S$ is a non-solvable subalgebra
of $\bar{L}$. This yields that $S=\bar{L}$ and hence $N=0$. Let
$A$ be a minimal ideal of $\bar{L}$. We have that either $A^2=A$
or $A^2=0$. In the latter case, we have $A\leq N=0$, a
contradiction. This yields that $A=A^2$ and hence $A=\bar{L}$.
Therefore $\bar{L}$ is simple. This completes the proof.

Next we consider the class of supersolvable Lie algebras

The structure of strongly-solvable minimal non-supersolvable Lie
algebras as well as that of $\phi$-free, non-strongly-solvable and
minimal non-supersolvable Lie algebras were determined in
\cite{ev}. For the readers' convenience, we give them below

\begin{theor}\cite[Theorems 1.1 and 1.2]{ev}\label{t:ev} Let $L$
be minimal non-supersolvable.
\begin{enumerate} \item If $L$ is strongly solvable,
then $L=\phi(L)\dot{+}A\dot{+}Fx$, where $A$ is subspace of $L$,
$A^2\leq \phi(L)$, with ${\rm ad}x$ acting irreducibly on $A$ and
${\rm ad}x\mid_{\phi(L)}$ is split.
\item If $L$ is $\phi$-free and non-strongly-solvable,
then ${\rm char}(F)=p>0$ and one of the following hold:
\begin{enumerate}
\item $L=((x,y,e_0,e_1,\ldots,e_{p-1}))$ with $[e_i,y]=(\alpha+i)e_i$
where $\alpha$ is a fixed scalar in $F$, $[e_i,x]=e_{i+1}$
(indices mod.$p$), $[x,y]=x$, $[e_i,e_j]=0$ and $F=\{t^p-t\mid
t\in F\}$.
\item $L=((x,y,z,e_0,e_1,\ldots,e_{p-1}))$ with $[e_i,z]=e_i$ for
every $i$ and $[e_i,x]=e_{i+1}$ $(i=0,\ldots,p-2)$ (non-specified
products are zero) and $F$ is perfect whenever $p=2$.
\end{enumerate}
\end{enumerate}
\end{theor}

Now we can prove the following

\begin{theor}\label{t:theor3} Let $L$ be a solvable Lie algebra.
\begin{enumerate} \item If $L$ minimal non-supersolvable, then $L$ is two-generated.
\item If every two-generated proper
subalgebra of $L$ is supersolvable but $L$ is not supersolvable,
then every proper subalgebra of $L$ is supersolvable. So, either
$L$ has the structure given in (1) of Theorem \ref{t:ev} or
$L/\phi(L)$ is isomorphic to one of the Lie algebras described in
(2) of Theorem \ref{t:ev}.
\end{enumerate}
\end{theor}

{\em Proof.} (1): If $L^2$ is not nilpotent, then $L$ is minimal
non-strongly-solvable. So, by Theorem \ref{t:theor1} it follows
that $L$ is two-generated. Now assume that $L^{2}$ is nilpotent.
Put $\bar{L}=L/\phi(L)$. We have that $\bar{L}$ is $\phi$-free,
strongly solvable and minimal non-supersolvable. By Theorem
\ref{t:ev}, we have that $\bar{L}=A\dot{+} Fx$ where $A$ is a
nonzero abelian ideal of $\bar{L}$ and ${\rm ad}\,x$ acts
irreducibly on $A$. Pick $0\neq a\in A$. We have that $\bar{L}$ is
generated by $a$ and $x$. From Lemma \ref{l:lemma2} it follows
that $L$ is two-generated. The proof of (1) is complete.

(2): Clearly, $L$ has a subalgebra $S$ which is minimal
non-supersolvable. From (1) it follows that $S$ is two-generated.
Then, by our hypothesis, we have that $S=L$.

\section{The property of triangulability}

A subalgebra $S$ of $L$ is said to be {\em triangulable} on $L$ if
${\rm ad}_{L} S = \{{\rm ad}_{L} x \mid x \in S\}$ is a Lie
algebra of linear transformations of $L$ which is triangulable
over the algebraic closure of $F$. A subalgebra $S$ of $L$ is said
to be {\em nil} on $L$ if ${\rm ad}\, x$ is nilpotent for every
$x\in L$. It is well-known that for every subalgebra $S$ of $L$,
there is a unique maximal ideal ${\rm nil}(S)$ of $S$ consisting
of ad-nilpotent elements of $L$ (see \cite[Proposition 2.1]{w}).
Also, it is known that $S$ is triangulable on $L$ if and only if
$S/{\rm nil}(S)$ is abelian (see \cite[Theorem 2.2]{w}). Note that
every subalgebra of $L$ which is triangulable on $L$ is strongly
solvable.

First, we give the following easy lemma which will be used in the
sequel

\begin{lemma}\label{l:nil} Let $L$ be a Lie algebra and $S$ and $T$ subalgebras of $L$
which are nil on $L$. Assume that $[S,T] \subset T$. Then $S + T$
is nil on $L$.
\end{lemma}

{\em Proof.} We see that the set ${\rm ad}\,S\cup{\rm ad}\,T$ is
weakly closed in the sense of Jacobson \cite{J}. Therefore, by
\cite[Theorem 1,p.33]{J}, every element of $S + T$ acts
nilpotently on $L$.

\begin{propo}\label{p:semi} Let $S$ be a triangulable subalgebra of
$L$. Then ${\rm nil}(S)$ is precisely the set of all elements of
$S$ which act nilpotently on $L$.
\end{propo}

{\em Proof.} Let $\Omega$ be an algebraic closure of $F$ and let
$L_{\Omega}=L\otimes_F\Omega$. Put $A={\rm ad}_LS$. We have that
$A\leq gl(L)$ and that $A_{\Omega}$ is a subalgebra of $
gl(L_{\Omega})$ and a set of simultaneously triangulable linear
maps. Then it is obvious that the set $N(A_{\Omega})$ of nilpotent
elements of $A_{\Omega}$ is closed under linear combinations and
that $(A_{\Omega})^2<N(A_{\Omega})$. Since $(A^2)_{\Omega}=
(A_{\Omega})^2$, it follows that $A^2$ is contained in the set $N$
of nilpotent elements in $A$. Also, we have that $N$ is closed
under linear combinations. So $N$ is actually an ideal of $A$.
Then ${\rm nil}(S)$ is contained in the inverse image $N^{\prime}$
of $N$ under the adjoint representation ${\rm ad}_L$. This is a
Lie-homomorphism, whence $N^{\prime}$ is an ideal of $S$. It
follows that $N^{\prime}={\rm nil}(S)$. This completes the proof.

\begin{lemma}\label{l:bige} Let $L$ be a Lie algebra which is not two-generated. Then
for each maximal subalgebra $M$ of $L$ and each non-zero element
$x$ in $M$ there exists a maximal subalgebra $S$ of $L$ different
from $M$ containing $x$. If the ground field $F$ is infinite, then
each element of $L$ lies in infinitely many maximal subalgebras of
$L$.
\end{lemma}

{\em Proof.} Let $M$ be a maximal subalgebra of $L$ and $0 \neq x
\in M$. Pick $y \in L \setminus M$. By hypothesis, $<x,y> \neq L$.
So, there exists a maximal subalgebra $S$ of $L$ containing
$<x,y>$. We see that $S \neq M$ since $y \in S$ and $y \not\in M$.
Now suppose that $F$ is infinite. Let $x \in L$ and assume that
$M_{1},\cdots, M_{r}$ are the only maximal subalgebras of $L$
which contain $x$. Since $F$ is infinite, $M_{1} \cup \cdots \cup
M_{r} \neq L$. Pick $y \in L$, $y \not\in M_{1} \cup \cdots \cup
M_{r}$. Then, by hypothesis $L \neq <x,y>$. Take a maximal
subalgebra $S$ of $L$ containing $<x,y>$. We see that $x \in S$
and $S \neq M_{i}$ for every $i$, which is a contradiction.

Now we give some information on the structure of a simple Lie
algebra all of whose proper subalgebras are triangulable. Clearly,
a simple minimal non-abelian Lie algebra satisfies this condition.
However, we do not know any other example.

\begin{theor}\label{t:simple-triang} Let $L$ be a simple Lie algebra such that every proper
subalgebra of $L$ is triangulable on $L$. Then the following hold:
\begin{enumerate}
\item if $M$ is a maximal subalgebra of $L$, then either $M$ is
abelian and none of the non-zero elements of $M$ acts nilpotently
on $L$ or ${\rm nil}(M)$ is a non-zero maximal nil subalgebra of
$L$;
\item if $K$ is a non-zero maximal nil subagebra of $L$, then the
normalizer of $K$ in $L$ is a maximal subalgebra of $L$;
\item for each two different maximal subalgebras $M_{1}$ and $M_{2}$
of $L$, none of the non-zero elements of $M_{1} \cap M_{2}$ acts
nilpotently on $L$ (in particular, $M_1\cap M_2$ is abelian); and
\item $L$ is two-generated.
\end{enumerate}
\end{theor}

{\em Proof.} (1): Let $M$ be any maximal subalgebra of $L$. First
assume that ${\rm nil}(M) = 0$. Then we have $M^2 \leq {\rm
nil}(M) = 0$. So, $M$ is abelian. By Proposition \ref{p:semi} it
follows that $M$ does not contain any non-zero ad-nilpotent
element of $L$. Now assume ${\rm nil}(M) \neq 0$. Let $K$ be a
maximal nil subalgebra of $L$ containing ${\rm nil}(M)$. Assume $K
\neq {\rm nil}(M)$. By Proposition \ref{p:semi} again, we have $M
\cap K = {\rm nil}(M)$. Nilpotency of $K$ implies that ${\rm
nil}(M)$ is properly contained in the normalizer $N_{K}({\rm
nil}(M))$ of ${\rm nil}(M)$ in $K$. This yields that ${\rm
nil}(M)$ is an ideal of $L$, since $L=M+N_{K}({\rm nil}(M))$ by
maximality of $M$, a contradiction. Therefore $K={\rm nil}(M)$ and
so ${\rm nil}(M)$ is maximal nil on $L$. The proof of (1) is
complete.

(2): Let $0 \neq K$ be a maximal nil subalgebra of $L$. Let $M$ be
a maximal subalgebra of $L$ containing $K$. Then, by Proposition
\ref{p:semi} it follows that $K = {\rm nil}(M)$. This yields that
$M =N_ {L}(K)$ and hence $N_{L}(K)$ is a maximal subalgebra.

(3): Let $M_1, M_2$ be distinct maximal subalgebras of $L$. By
Proposition \ref{p:semi} we have that $${\rm nil}(M_1\cap
M_2)={\rm nil}(M_1)\cap {\rm nil}(M_2)$$We need to prove that
${\rm nil}(M_1\cap M_2)=0$. Assume that ${\rm nil}(M_1\cap
M_2)\neq 0$. Choose $M_1$ and $M_2$ such that the dimension of
${\rm nil}(M_1\cap M_2)$ is maximal. Suppose that ${\rm
nil}(M_1\cap M_2)={\rm nil}(M_2)$. Then we have $0\neq {\rm
nil}(M_2)\leq {\rm nil}(M_1)$. By (1) it follows that ${\rm
nil}(M_2)={\rm nil}(M_1)$. This yields that ${\rm nil}(M_2)$ is an
ideal of $M_1$ and $M_2$ and therefore it is an ideal of $L$.
Simplicity of $L$ implies that ${\rm nil}(M_2)=0$, a
contradiction. Hence ${\rm nil}(M_1\cap M_2)\neq {\rm nil}(M_2)$.
Analogously, we have ${\rm nil}(M_1\cap M_2)\neq {\rm nil}(M_1)$.
Now, by using Engel's Theorem we obtain subalgebras $S_1$ and
$S_2$ of $L$ strictly containing ${\rm nil}(M_1\cap M_2)$ and such
that
$${\rm nil}(M_1\cap M_2)\lhd S_i\leq{\rm nil}(M_i)\,\,{\rm
for}\,\, i=1,2$$ Let $S$ be the subalgebra of $L$ generated by
$S_1$ and $S_2$. We find that ${\rm nil}(M_1\cap M_2)$ is an ideal
of $S$. Simplicity of $L$ implies that $S\neq L$. Take a maximal
subalgebra $M_3$ of $L$ containing $S$. We have $${\rm
nil}(M_1\cap M_2)\subset S_1\subseteq M_3\cap{\rm nil}(M_1)={\rm
nil}(M_1\cap M_3)$$ which contradicts the maximality of ${\rm
nil}(M_1\cap M_2)$. This completes the proof of (3).

(4): Assume that $L \neq <a,b>$ for every $a$, $b \in L$. Let $K$
be a nil subalgebra of $L$ of maximal dimension. If $K = 0$, then
by (1) we have that every maximal subalgebra of $L$ is abelian.
Let $M_1,M_2$ be distinct maximal subalgebras of $L$. We have that
$M_1\cap M_2=0$; since, otherwise, we would have that $M_1\cap
M_2$ is a nonzero ideal of $L$, which is a contradiction. Pick
$0\neq a_1\in M_1$ and $0\neq a_2\in M_2$. Since $ <a_1,a_2>\neq
L$, we can take a maximal subalgebra $M_3$ of $L$ containing $a_1$
and $a_2$. But then we have $0\neq a_1\in M_1\cap M_3$ and
$M_1\neq M_3$, a contradiction. Therefore $K \neq 0$. Pick $x \in
K$, $x \neq 0$. By (3), there exists only one maximal subalgebra
of $L$ containing $x$, which contradicts Lemma \ref{l:bige}.
Therefore $L$ can be generated by two elements. The proof of the
theorem is complete.

Now we can prove the following

\begin{lemma}\label{l:t}
\begin{enumerate}
\item If $\phi(L)\leq S\leq L$, then
${\rm nil}(S/\phi(L))={\rm nil}(S)/\phi(L)$.
\item If every two-generated {\em proper} subalgebra of $L$ is
triangulable on $L$, then the same holds in $L/\phi(L)$.
\end{enumerate}
\end{lemma}

{\em Proof.} (1): Assume $\phi(L)\leq S\leq L$. Let ${\rm
nil}(S/\phi(L))=K/\phi(L)$. Since $\phi(L)$ is a nilpotent ideal
of $L$ (see \cite{to}), it is nil on $L$. So, $\phi(L)\leq {\rm
nil}(S)$. Clearly, ${\rm nil}(S)\leq K$. Now let $x\in K$. We have
that $x$ acts nilpotently on $L/\phi(L)$. This yields that
$L=\phi(L)+L_0(x)$, where $L_0(x)$ is the Fitting null-component
of $L$ relative to ${\rm ad}\,x$. As $L_0(x)$ is a subalgebra of
$L$, it follows that $L=L_0(x)$. So that $x$ acts nilpotently on
$L$. This yields that $K\leq{\rm nil}(S)$.

(2): Assume that every two-generated {\em proper} subalgebra of
$L$ is triangulable on $L$. Let $S/\phi(L)$ be a two-generated
proper subalgebra of $L/\phi(L)$. Take elements $a,b\in L$ such
that $S/\phi(L)=<a+\phi(L),b+\phi(L)>$. We have that
$S=<a,b>+\phi(L)$ and that $<a,b>$ is triangulable on $L$. On the
other hand, from Lemma \ref{l:nil} it follows that $\phi(L)+{\rm
nil}(<a,b>)$ is nil on $L$. Since $\phi(L)+{\rm nil}(<a,b>)$ is an
ideal of $S$, we have that $\phi(L)+{\rm nil}(<a,b>)\leq {\rm
nil}(S)$. Then,
$$S^2\leq \phi(L)+<a,b>^2\leq \phi(L)+{\rm nil}(<a,b>)\leq
{\rm nil}(S)$$ So, we have that
$$(S/\phi(L))^2=S^2+\phi(L)/\phi(L)\leq{\rm nil}(S)/\phi(L)={\rm nil}(S/\phi(L))$$
by (1). Therefore $S/\phi(L)$ is triangulable on $L/\phi(L)$. The
proof is complete.

\begin{theor}\label{t:2}
\begin{enumerate}
\item If $L$ is solvable and every two-generated proper
subalgebra of $L$ is triangulable on $L$, then $L$ is
triangulable.
\item If every proper subalgebra of $L$ is triangulable on $L$ but $L$ itself is not,
then $L$ is two-generated and $L/{\phi(L)}$ is simple (so,
$L/\phi(L)$ satisfies (1)-(4) in Theorem \ref{t:simple-triang}).
\end{enumerate}
\end{theor}

{\em Proof.} (1): Let $L$ be solvable and every two-generated
proper subalgebra of $L$ be triangulable on $L$. In particular, we
have that every two-generated proper subalgebra of $L$ is strongly
solvable. Put $\bar{L}=L/\phi(L)$. Then, by Theorem
\ref{t:theor1}(1) it follows that every proper subalgebra of
$\bar{L}$ is strongly solvable. Moreover, by Lemma \ref{l:t} we
have that every two-generated proper subalgebra of $\bar{L}$ is
triangulable on $\bar{L}$. Now suppose that $L$ is not
triangulable. So that $L$ is not strongly solvable. Thus,
$\bar{L}$ is not strongly solvable either. By using Theorem
\ref{t:bt} we obtain that $\bar{L}=A+ B$, where $A$ is the unique
minimal ideal of $\bar{L}$, $A^2=0$, $B<\bar{L}$, $A\cap B=0$, and
either $B=M\dot{+} Fx$, where $M$ is the unique minimal ideal of
$B$ or $B$ is the three-dimensional Heisenberg algebra. We have
that $\bar{L}^2=A+B^2$. On the other hand, we see that $B$ is
two-generated. So, $B$ is triangulable on $\bar{L}$. This yields
that $B^2$ acts nilpotently on $A$ and therefore $\bar{L}^2$ is
nilpotent. So that $\bar{L}$ is strongly solvable. This
contradiction completes the proof of (1).

(2): Assume that every proper subalgebra of $L$ is triangulable on
$L$ but $L$ itself is not. By (1) we have that $L$ is not
solvable. By Lemmas \ref{l:t}(1), \ref{l:lemma2} and since the
class of solvable Lie algebras is saturated, we may suppose
without loss of generality that $L$ is $\phi$-free. By Towers
\cite{to} we have that $L=N+S$ and $N\cap S=0$, where $N$ is the
largest nilpotent ideal of $L$ and $S$ is a subalgebra of $L$. If
$S\neq L$, then we have that $S$ is solvable. So, $L$ is solvable,
a contradiction. It follows that $N=0$ and so $L$ is simple.
Hence, $L$ satisfies (1)-(4) in Theorem \ref{t:simple-triang}. The
proof is complete.
\bigskip

ACKNOWLEDGMENT

The authors are grateful to the referee for his/her useful
comments; in particular, for providing a shorter and clearer
proof of Proposition 4.2.

\end{document}